\pdfoutput=1  
\documentclass[runningheads]{llncs}

\usepackage[utf8]{inputenc}
\usepackage[T1]{fontenc}
\usepackage{lmodern}
\usepackage{amsmath,amssymb}
\usepackage{mathtools}
\usepackage{graphicx}
\usepackage{booktabs}
\usepackage{array}
\usepackage{caption}
\usepackage{subcaption}
\usepackage{multirow}
\usepackage{float}
\usepackage{enumitem}
\usepackage[hidelinks]{hyperref}
\usepackage{xcolor}

\hypersetup{
  pdftitle={A General Variable Neighborhood Search for Logistic Network Flow
            Design with Nonlinear Edge Costs and Transshipment Capacities},
  pdfauthor={Varvara Kovaleva, Alexander Ponomarenko},
  pdfsubject={Multi-commodity network flow design; variable neighborhood search},
  pdfkeywords={multi-commodity flow, network design, variable neighborhood
               search, metaheuristics, mixed integer programming, logistics}
}

\usepackage{placeins}  
\usepackage[ruled,vlined,linesnumbered]{algorithm2e}
\SetKwProg{Fn}{Function}{:}{}
\SetKwFunction{PN}{ProcessNeighborhood}
\SetKwFunction{GenPaths}{GenPaths}
\SetKwFunction{SolveMILP}{SolveMILP}
\SetKwFunction{RoundedLP}{RoundedLP}

\newcommand{\E}{E}
\newcommand{\Vv}{V}
\newcommand{\Real}{\mathbb{R}}
\newcommand{\Int}{\mathbb{Z}}

\begin{document}

\title{A General Variable Neighborhood Search for Logistic Network
Flow Design with Nonlinear Edge Costs and Transshipment Capacities}
\titlerunning{GVNS for Logistic Network Flow Design}

\author{Varvara Kovaleva \inst{1} \and Alexander Ponomarenko\inst{2} }
\authorrunning{V. Kovaleva and A. Ponomarenko}
\institute{
HSE University, Nizhny Novgorod, Russia\\
\email{varvara\_kovaleva03@mail.ru}
\and
Laboratory of Algorithms and Technologies for Network Analysis, HSE University, Nizhny Novgorod, Russia\\
\email{aponomarenko@hse.ru}\\
}

\maketitle

\begin{abstract}
We study the problem of routing multiple commodities through a logistic
transport network so as to minimise the joint cost of vehicle movements and
of transshipment (reload) operations at intermediate storages. The problem
combines two features that are common in practice but rarely treated
together: edge costs that are \emph{nonlinear} in the transported volume,
because freight is carried by vehicles of fixed capacity, and \emph{hard
upper bounds} on the volume that can be reloaded at each node. We formulate
the problem as a mixed integer linear program in both an arc-based and a
path-based form, derive its linear relaxation, and propose a metaheuristic
solver based on General Variable Neighborhood Search (GVNS). The method works
in the path space and relies on three problem-specific neighborhoods that
generate new candidate paths by (i)~shortest paths under a reload-aware
metric, (ii)~edges carrying partially loaded vehicles, and (iii)~detours
around congested storages. We evaluate the method on our open
benchmark~\cite{benchmark} for this problem: a corpus of realistic logistic
networks calibrated on real road graphs of nine European countries, five U.S.\
states and the European part of Russia, with near-planar topology
(Zipf city sizes, a spatial-network budget model, and a doubly constrained
gravity demand model). On 35 synthetic and 15 real-geometry instances with up
to 150 storages, GVNS returns solutions within about one percent of the best
dual bound on small networks---a branch-and-bound lower bound substantially
tighter than the LP relaxation---is stable across random seeds, and for a
comparable wall-clock budget is up to $55\%$ better than a time-limited exact
solver on the medium instances.

\end{abstract}

\section{Introduction}

Organising commodity flows between distribution centres is one of the central
problems of logistic networks. The growth in scale and complexity of modern
warehousing structures creates a persistent demand for efficient algorithms
that build optimal transportation plans. Two characteristics make realistic
instances hard. First, the cost of using a road is \emph{nonlinear} in the
volume transported along it: goods travel in vehicles of fixed capacity, so
the relevant quantity is the integer \emph{number of vehicles}, obtained by
rounding the carried volume up to the capacity. Second, every storage has a
limited \emph{transshipment capacity}---an upper bound on the total volume
that may be unloaded, stored and reloaded there. Together these features turn
an otherwise classical multi-commodity flow problem into a mixed integer
program whose feasible region is shaped by node capacities and whose
objective contains a step-wise vehicle term.

Exact methods for multi-commodity flow problems are computationally expensive
because the problem is NP-hard. Column generation~\cite{barnhart1998} is a
classical approach that solves path-based formulations by dynamically building
a set of promising paths, and is frequently combined with branch-and-price. In
practice, heuristics and metaheuristics---greedy algorithms, genetic
algorithms, tabu search, simulated annealing and local-search-based
methods---are widely used to obtain good solutions in acceptable
time~\cite{crainic2015}. Among them, Variable Neighborhood Search (VNS) and
its variants have proven effective across a broad range of optimisation
problems. Masri et al.~\cite{masri2015} applied a multi-start VNS to the
single-path multi-commodity flow problem, and General VNS has performed well
on the capacitated vehicle routing problem and other flow-constrained logistic
problems~\cite{amous2017}. Closer to the present setting, problem-specific
heuristics have proven effective on real-life logistic routing problems, such as
the truck-and-trailer routing problem with time windows and split
deliveries~\cite{batsyn2014ttrp,batsyn2015heuristic}.

\paragraph{Contributions.} The aim of this work is the development and study
of a metaheuristic algorithm, based on General Variable Neighborhood Search
(GVNS), for distributing commodity flows in a logistic network with nonlinear
edge costs and limited storage capacities. Specifically, we
\begin{enumerate}[leftmargin=1.4em,itemsep=1pt,topsep=2pt]
\item formalise the problem and give equivalent arc-based and path-based MILP
      models together with their linear relaxation (Section~\ref{sec:model});
\item design a GVNS that searches in the space of commodity paths and is
      equipped with three problem-specific neighborhoods that exploit reload
      costs, partial vehicle loads and storage congestion
      (Section~\ref{sec:gvns});
\item report a computational study on our open benchmark~\cite{benchmark}
      (Section~\ref{sec:data}), comparing GVNS
      against the LP relaxation, its rounded version, and MIP dual bounds and
      incumbents from an exact solver
      (Section~\ref{sec:experiments}).
\end{enumerate}
The main methodological elements are the specialised GVNS neighborhoods adapted
to the path-based formulation and its constraints, together with the
exact-solver study that quantifies the quality of the resulting solutions.

\section{Problem statement and model}
\label{sec:model}

\subsection{Verbal statement}
We are given a transport network of storages connected by roads. Prescribed
volumes of goods must be moved from origins to destinations so as to minimise
the total cost of moving vehicles along roads plus the cost of transshipment
(unloading, storing and reloading) at intermediate transit storages. Each
storage has a limited transshipment capacity and a per-unit reload cost.
Freight travels in vehicles of fixed capacity, and each road has a cost per
vehicle pass; the number of vehicles on a road must be sufficient to carry the
total volume routed through it. The problem belongs to the class of
multi-commodity flow problems with additional node-capacity constraints and a
vehicle-induced nonlinear edge cost, and admits two equivalent formulations:
arc-based and path-based.

\subsection{Data}
The road network is a directed graph $G=(\Vv,\E)$, where $\Vv$ is the set of
storages and $\E$ the set of roads. We are given a set of transportation
requests (commodities) $K$; for each $k\in K$ we know an origin $s_k\in \Vv$,
a destination $t_k\in \Vv$, and a demand $d_k\in\Real_+$ to be shipped from
$s_k$ to $t_k$. For each storage $v\in \Vv$ we know the per-unit reload cost
$f_v\in\Real_+$ and the maximum admissible transshipment volume
$W_v\in\Real_+$. Each edge $e\in \E$ has a per-vehicle traversal cost
$c_e\in\Real_+$. Finally, $C\in\Real_+$ is the capacity of one vehicle.

\subsection{Decision variables}
\textit{Arc-based formulation.} Let $x^k_e\in\Real_+$ be the flow of commodity
$k\in K$ on edge $e\in \E$ and $y_e\in\Int_+$ the number of vehicles traversing
$e$. The model has $|K|\cdot|\E|+|\E|$ variables.

\textit{Path-based formulation.} Let $\Pi_k$ be the set of all paths from $s_k$
to $t_k$. Let $x^k_p\in\Real_+$ be the flow of commodity $k$ on path
$p\in\Pi_k$, and $y_e\in\Int_+$ as above. The model has
$\sum_{k\in K}|\Pi_k|+|\E|$ variables.

\subsection{Objective}
Denote by $\delta^-(v)\subseteq \E$ and $\delta^+(v)\subseteq \E$ the sets of
edges entering and leaving $v$. The objective minimises the sum of vehicle
movement cost (first term) and transit reload cost (second term).

\noindent\emph{Arc-based:}
\begin{equation}
\min\;\Bigg(\sum_{e\in \E} c_e\,y_e
 \;+\!\!\sum_{k\in K}\;\sum_{\substack{v\in \Vv\\ v\neq s_k,\,v\neq t_k}}
 \sum_{e\in\delta^-(v)} f_v\,x^k_e\Bigg).
\label{eq:obj-arc}
\end{equation}

\noindent\emph{Path-based:}
\begin{equation}
\min\;\Bigg(\sum_{e\in \E} c_e\,y_e
 \;+\!\!\sum_{k\in K}\;\sum_{p\in\Pi_k}\;
 \sum_{\substack{v\in p\\ v\neq s_k,\,v\neq t_k}} f_v\,x^k_p\Bigg).
\label{eq:obj-path}
\end{equation}

\subsection{Constraints}
\textbf{(1) Demand satisfaction.}
\begin{align}
\text{arc:}\quad
&\sum_{e\in\delta^+(v)} x^k_e-\sum_{e\in\delta^-(v)} x^k_e=
 \begin{cases} d_k, & v=s_k\\ -d_k, & v=t_k\\ 0, & \text{otherwise}\end{cases}
 \quad \forall k\in K,\ \forall v\in \Vv;
 \label{eq:demand-arc}\\[2pt]
\text{path:}\quad
&\sum_{p\in\Pi_k} x^k_p=d_k, \qquad \forall k\in K.
\label{eq:demand-path}
\end{align}

\textbf{(2) Vehicle capacity on edges.} Enough vehicles must run on each edge
to carry the total flow:
\begin{align}
\text{arc:}\quad &\sum_{k\in K} x^k_e\le C\,y_e, \qquad \forall e\in \E;
\label{eq:cap-arc}\\
\text{path:}\quad &\sum_{k\in K}\sum_{p\in\Pi_k:\,e\in p} x^k_p\le C\,y_e,
\qquad \forall e\in \E.
\label{eq:cap-path}
\end{align}

\textbf{(3) Storage transshipment capacity.}
\begin{align}
\text{arc:}\quad
&\sum_{\substack{k\in K\\ v\neq s_k,\,v\neq t_k}}\sum_{e\in\delta^-(v)} x^k_e
 \le W_v, \qquad \forall v\in \Vv;
\label{eq:node-arc}\\
\text{path:}\quad
&\sum_{\substack{k\in K\\ v\neq s_k,\,v\neq t_k}}\sum_{p\in\Pi_k}
 \sum_{v\in p} x^k_p \le W_v, \qquad \forall v\in \Vv.
\label{eq:node-path}
\end{align}

\subsection{Classification and linear relaxation}
\label{sec:relax}
The objective and all constraints are linear, and the variables are mixed:
the vehicle counts $y_e$ are integer while the flows $x^k_e$ (or $x^k_p$) are
continuous. The problem is therefore a \emph{mixed integer linear program}
(MILP), an NP-hard class for which no polynomial-time exact algorithm is known
in general.

A standard way to bound the optimum is the LP relaxation: dropping the
integrality of $y_e$ turns the MILP into a linear program solvable in
polynomial time, and its optimum is a lower bound (for minimisation). Once
$y_e$ is continuous, the capacity constraints \eqref{eq:cap-arc}--%
\eqref{eq:cap-path} hold with equality at optimum, so $y_e$ can be eliminated:
\begin{equation}
y_e=\frac1C\sum_{k\in K} x^k_e \quad(\text{arc}), \qquad
y_e=\frac1C\sum_{k\in K}\sum_{p\in\Pi_k:\,e\in p} x^k_p \quad(\text{path}).
\label{eq:ye-relax}
\end{equation}
Substituting \eqref{eq:ye-relax} into the objective yields a relaxed problem
in the flow variables only. In the path-based form,
\begin{equation}
\min\;\sum_{k\in K}\sum_{p\in\Pi_k} x^k_p
\Bigg(\frac1C\sum_{e\in p} c_e+\sum_{\substack{v\in p\\ v\neq s_k,\,v\neq t_k}}
f_v\Bigg),
\label{eq:relax-path}
\end{equation}
which makes explicit that an optimal relaxed solution prefers paths with low
travel cost and few transit nodes.

\section{The GVNS solution approach}
\label{sec:gvns}

\subsection{General Variable Neighborhood Search}
VNS~\cite{hansen2019,hansen2017} is a metaheuristic whose core idea is the
systematic change of neighborhood structures while exploring the solution space.
Given a problem $\min\{f(x):x\in X\}$ and a family of neighborhoods
$N_k(x),\ k=1,\dots,k_{\max}$, the classical scheme alternates three steps: a
\emph{shake} (a random move in $N_k$ to escape local optima), a \emph{local
search}, and a \emph{neighborhood change} (on improvement, accept the move and
reset $k\!\gets\!1$; otherwise advance to $k\!+\!1$). When the local search is
itself a Variable Neighborhood Descent (VND)---a deterministic sweep through the
neighborhoods ordered from intensification to diversification---the method is
called General VNS (GVNS), which often outperforms single-neighborhood
heuristics~\cite{karakostas2019,mjirda2014}. Our instantiation searches the space
of commodity paths; its concrete form is Algorithm~\ref{alg:gvns}, built on the
neighborhood-processing primitive of Algorithm~\ref{alg:pn}.

\subsection{Preprocessing and initial solution}
Transport networks are typically sparse directed graphs in which not every
pair of storages is joined by a direct road, and not every transit node needs
a reload. Working on the raw graph both shrinks the feasible region (through
node capacities) and may even render the instance infeasible if no path exists
between an origin and a destination.

To avoid this, the graph is \emph{completed}: for every pair of storages
$v,u$ we add (or update) an edge whose cost equals that of the shortest path
from $v$ to $u$ in $G$. Each edge of the resulting complete graph thus
represents the cheapest stop-free route between its endpoints. This guarantees
the existence of a feasible solution, because flow routed on the direct
origin--destination edge does not pass through any transit node and hence
violates no node-capacity constraint. From the relaxed objective
\eqref{eq:relax-path} it follows that the cheapest, transit-free direct edges
are exactly the ones a good solution prefers. We therefore take as the initial
solution the relaxed LP solution that sends each commodity entirely along its
direct origin--destination edge, with vehicle counts rounded up:
\begin{equation}
y_e=\Big\lceil \tfrac1C\sum_{k\in K} x^k_e\Big\rceil, \qquad \forall e\in \E.
\label{eq:init}
\end{equation}

\subsection{Neighborhoods}
Since the vehicle counts $y_e$ are determined directly from the flows via
\eqref{eq:init}, the neighborhoods vary the flows. Predicting feasible
per-edge flows is hard because edges are shared by many commodities and
adjacent flows are strongly coupled, so the whole network must be considered
at once. A single path, by contrast, describes the routing of just one
commodity, and adding a new path to an already feasible solution never worsens
it. The neighborhoods therefore \emph{build and modify the set of candidate
paths for each commodity}: for every commodity a set of new paths is generated
and added to the model, a MILP solver is run for a limited time to
redistribute flow, and---if the objective improves---all new paths carrying
nonzero flow are kept. Three neighborhoods define how the new paths are
generated.

\paragraph{N1 -- Shortest paths (reload-aware).}
Classical shortest-path search ignores reload costs. We instead build an
auxiliary complete graph $G'$ whose edge weights measure the per-unit cost of
moving one unit of commodity, including reloads at the endpoints:
\begin{equation}
c'_{vu}=f_v+f_u+\frac{c_{vu}}{C}, \qquad \forall (v,u)\in \E .
\end{equation}
Before GVNS starts, a pool of shortest paths under $c'$ is precomputed for
each commodity; new paths are then drawn from this pool, either at random or
preferring low-cost paths. The neighborhood parameters are the number of new
paths per commodity and a boolean choosing the selection strategy. During VND
all combinations are swept (one to three new paths, both strategies); in Shake,
five random paths are always added, which helps escape local optima.

\paragraph{N2 -- Partially loaded edges.}
A partially loaded vehicle has spare room: routing extra flow over such an edge
can free another vehicle and reduce the vehicle count. For every edge we
compute the residual load
\begin{equation}
load_e=\Big(\sum_{k\in K}\sum_{p\in\Pi_k:\,e\in p} x^k_p\Big) \bmod C,
\qquad \forall e\in \E,
\end{equation}
and build a graph $G''$ that favours edges with free capacity. Fully loaded
edges ($load_e=0$) are least attractive, since new flow would require an extra
vehicle, and keep their weight from $G'$. The remaining edges receive
\begin{equation}
c''_{vu}=f_v+f_u+\frac{load_{vu}}{C}\Big(\frac{c_{vu}}{C}\Big)^{\!\alpha},
\qquad \alpha\le 1,
\end{equation}
so that smaller residual load (more free room) lowers the weight and raises
priority; $\alpha$ tunes how strongly travel cost influences the weight.
Shortest-path search on $G''$ then tends to build paths through
partially loaded vehicles. The number of new paths is set dynamically: a
counter of fully loaded paths and an accumulator of total free room stop
generation once a threshold is reached (in VND, at most two fully loaded paths
or total free room exceeding one vehicle capacity; in Shake, three paths and
$3C$).

\paragraph{N3 -- Detours around congested storages.}
Storages whose transshipment volume is at or near their capacity are
bottlenecks, and routing around them is desirable. For the current solution we
compute a congestion ratio
\begin{equation}
overload_v=\frac1{W_v}\sum_{\substack{k\in K\\ v\neq s_k,\,v\neq t_k}}
\sum_{p\in\Pi_k}\sum_{v\in p} x^k_p, \qquad \forall v\in \Vv .
\end{equation}
Storages are sorted by decreasing $overload_v$, and a subset of the most
congested ones is selected according to a threshold (and trimmed or padded to a
target size). For each prefix $\{v_1,\dots,v_i\}$ of the ordered set, those
nodes (except the commodity's own origin and destination) are deleted from
$G'$ and a shortest path is recomputed, yielding detours that bypass the
busiest nodes and shift flow to less loaded ones. In VND the subset size ranges
from $|\Vv|/5$ to $|\Vv|/4$, with threshold $0.75$ and one path per
configuration; in Shake the range widens to $|\Vv|/4$--$|\Vv|/3$, the threshold
drops to $0.5$ and three paths are generated per configuration.

\begin{algorithm}[t]
\DontPrintSemicolon
\Fn{\PN{$P,\ x,\ z,\ N,\ \mathit{div},\ \tau$}}{
  \tcp{incumbent path pool $P$, flows $x$, objective $z=f(x)$; updated in place}
  $Q \gets \bigcup_{k \in K}{}$ \GenPaths{$N,\ k,\ \mathit{div}$}\tcp*{candidate paths; $N$ selects the rule}
  \lIf{$Q = \emptyset$}{\Return \textbf{false}}
  $(x',\ z') \gets{}$ \SolveMILP{$P \cup Q,\ x,\ \tau$}\tcp*{redistribute flow; warm start $x$, limit $\tau$}
  \eIf{$z' < z$}{
    $P \gets P \cup \{\, p \in Q : x'_p > 0 \,\}$;\ \ $x \gets x'$;\ \ $z \gets z'$\;
    prune from $P$ every path carrying zero flow\tcp*{keep the model bounded}
    \Return \textbf{true}
  }{
    discard $Q$\;
    \Return \textbf{false}
  }
}
\caption{\texttt{ProcessNeighborhood}: enlarge the path pool, re-optimise the
flows by MILP, and keep the move iff it improves the objective.}
\label{alg:pn}
\end{algorithm}

\begin{algorithm}[t]
\DontPrintSemicolon
\SetKwInOut{Input}{Input}\SetKwInOut{Output}{Output}
\Input{graph $G=(\Vv,\E)$, commodities $K$, vehicle capacity $C$, time limit $t_{\max}$}
\Output{path pool $P$ with commodity flows $x$}
build the reload-aware complete graph $G'$; precompute a shortest-path pool per commodity\;
$P \gets \{\,(s_k,t_k) : k \in K\,\}$;\ \ $(x,z) \gets{}$ \RoundedLP{$P$}\tcp*{initial solution, Eq.~\eqref{eq:init}}
$k \gets 1$\;
\While{$k \le 3$ \textbf{\emph{and}} $\mathrm{CpuTime} < t_{\max}$}{
  \PN{$P,\ x,\ z,\ N_k,\ \textbf{true},\ \tau_{\mathrm{shk}}$}\tcp*{shake: diversified, larger budget}
  $z_0 \gets z$\;
  \Repeat{$\neg\,\mathit{imp}$ \textbf{\emph{or}} $\mathrm{CpuTime} \ge t_{\max}$}{
    $\mathit{imp} \gets{}$ \PN{$P,\ x,\ z,\ N_1,\ \textbf{false},\ \tau_{\mathrm{vnd}}$}\;
    \lIf{$\neg\,\mathit{imp}$}{$\mathit{imp} \gets{}$ \PN{$P,\ x,\ z,\ N_2,\ \textbf{false},\ \tau_{\mathrm{vnd}}$}}
    \lIf{$\neg\,\mathit{imp}$}{$\mathit{imp} \gets{}$ \PN{$P,\ x,\ z,\ N_3,\ \textbf{false},\ \tau_{\mathrm{vnd}}$}}
  }
  \leIf{$z < z_0$}{$k \gets 1$}{$k \gets k+1$}\tcp*{neighborhood change}
}
\Return $(P,\ x)$
\caption{\texttt{GVNS} for logistic network flow design.}
\label{alg:gvns}
\end{algorithm}

\subsection{Overall algorithm}
The method is given by Algorithms~\ref{alg:pn} and~\ref{alg:gvns}.
Algorithm~\ref{alg:pn} is the primitive shared by the three neighborhoods: it
generates a set of candidate paths---the only step that differs between $N_1$,
$N_2$ and $N_3$, via \GenPaths---enlarges the path pool, re-optimises the
commodity flows with the path-based MILP of Section~\ref{sec:model}, warm-started
from the current solution, and, \emph{only} if the objective improves, commits
the new paths that carry flow while pruning those that do not, which keeps the
model bounded. Algorithm~\ref{alg:gvns} embeds this primitive in the GVNS
template: starting from the rounded-LP solution of Eq.~\eqref{eq:init}, each round
shakes with the current neighborhood---diversified parameters and a larger MILP
budget $\tau_{\mathrm{shk}}$---and then runs the VND sweep
$N_1\!\to\!N_2\!\to\!N_3$ with a smaller budget $\tau_{\mathrm{vnd}}$, since VND
explores a narrower region. The neighborhood index resets on improvement and
advances otherwise; a global limit $t_{\max}$ bounds the run, and the search
stops once a full sweep brings no improvement.

\section{Benchmark instances}
\label{sec:data}

We evaluate GVNS on our open benchmark~\cite{benchmark,benchmark-data} for this
problem, published as a dataset on Zenodo (CC-BY-4.0); we summarise here only what
is needed for the study and refer the reader to~\cite{benchmark} for the full
generator. The corpus targets upper-level road networks connecting regional
sorting centres, with a generator calibrated on OpenStreetMap~\cite{osm} data for
nine European countries, five large U.S.\ states and the European part of Russia.
It reproduces the structural statistics of real regional transport
systems---Zipf-distributed city sizes~\cite{gabaix1999}, a spatial-network budget
topology~\cite{gastner2006}, centrality-driven transshipment
capacities~\cite{barthelemy2011,kirkley2018} and a doubly constrained gravity
demand model~\cite{wilson1967}---yielding \emph{near-planar, large-diameter}
graphs (mean vertex degree about $3.575$), the opposite of small-world networks.
The corpus comprises $35$ synthetic instances ($10$ with $10$ storages, $7$ each
with $20$ and $30$, $5$ with $50$, and $3$ each with $100$ and $150$) and $15$
real-geometry instances built on the calibration road graphs ($11$--$90$
storages); we use all $50$, each exactly reproducible from its published seed and
parameters.

\section{Computational experiments}
\label{sec:experiments}

\subsection{Setup, baselines and metric}
All experiments were run on the HPC cluster ``cHARISMa''~\cite{charisma} (nodes with two
AMD~EPYC CPUs, 44~cores and $\ge$720\,GB RAM per node, Rocky~Linux~9). Each job
was pinned to eight cores. The exact models were solved with
HiGHS~1.14.0~\cite{highs} and, as an independent cross-check of the lower
bounds, with SCIP~10.0.2~\cite{scip}; the metaheuristic uses HiGHS as its
in-loop redistribution solver, as before. GVNS is stochastic (random path
choices in the shake step), so every instance was solved with three independent
random seeds and we report the best objective together with its across-seed
mean and standard deviation. Unless stated otherwise the exact solves were run
to a target MIP gap of zero (i.e.\ to optimality where possible) under
a size-dependent wall-clock budget (30~min for $|V|=10$, growing to 4~h for
$|V|=150$); GVNS is run with a time limit of twice the number of storages, in
minutes. The two budgets are comparable---at $|V|=50$, for example, GVNS is
allowed $100$~min against the exact solver's $2$~h, and at $|V|=10$ it is
$20$~min against $30$~min.

The main algorithm settings are collected in Table~\ref{tab:params}. The size of
the arc-based MIP grows steeply with $|V|$---from about $6\,000$ variables at
$|V|=10$ (edges $|E|=90$, commodities $|K|=67$) to $3.7\times10^{8}$ at $|V|=150$
($|E|=22\,350$, $|K|=16\,762$)---which is why the exact model becomes intractable
well before that size while GVNS, working in the much smaller path space, remains
applicable throughout. A GVNS run took, on average, from about $640$~s at
$|V|=10$ to $3.7\times10^{3}$~s at $|V|=50$ and $1.8\times10^{4}$~s at $|V|=150$.

\begin{table}[t]
\centering
\caption{Main GVNS and exact-solver settings. VND/Shake values are given as
VND\,/\,Shake where they differ.}
\label{tab:params}
\small
\begin{tabular}{ll}
\toprule
Setting & Value\\
\midrule
Neighborhoods $k_{\max}$ & $3$ ($N_1,N_2,N_3$)\\
GVNS time limit & $2|V|$ minutes\\
Random seeds & $3$\\
Vehicle capacity $C$ & $50$\\
$N_1$ new paths per commodity (VND\,/\,Shake) & $1$--$3$\,/\,$5$\\
$N_2$ travel-cost exponent $\alpha$ & $0.5$\\
$N_2$ free-room stop (VND\,/\,Shake) & $1C$\,/\,$3C$\\
$N_3$ congestion threshold (VND\,/\,Shake) & $0.75$\,/\,$0.5$\\
$N_3$ storage-subset size (VND\,/\,Shake) & $|V|/5$--$|V|/4$\,/\,$|V|/4$--$|V|/3$\\
Exact-solver MIP-gap target & $0$\\
Exact-solver budget & $30$~min ($|V|{=}10$) \dots\ $4$~h ($|V|{=}150$)\\
\bottomrule
\end{tabular}
\end{table}

Lacking known optima, we compare GVNS against four references:
\begin{itemize}[leftmargin=1.4em,itemsep=1pt,topsep=2pt]
\item \textbf{LP relaxation} -- the optimum of the relaxed problem
      (Section~\ref{sec:relax}); a \emph{lower bound}. On the completed graph it
      is attained by routing every commodity along its direct edge
      (Section~\ref{sec:gvns}), so it has a closed form and is available at every
      size, though---ignoring vehicle integrality---it is loose;
\item \textbf{Rounded LP} -- the relaxed solution with $y_e$ rounded up; this is
      the GVNS starting solution and a valid \emph{upper bound};
\item \textbf{MIP dual bound} -- the best lower bound proved by the exact
      branch-and-bound solver within its budget, i.e.\ the global lower bound of
      the search tree: the LP relaxation strengthened by branching and cutting
      planes, which rises towards the optimum as the tree is explored. It is thus
      a \emph{valid} lower bound ($\underline{z}\le z^{\ast}$) that---unlike the
      root LP relaxation, which ignores vehicle integrality---is typically tight
      (within a few percent of the optimum); it is the key reference, letting us
      bound the true optimum even when optimality is not proved;
\item \textbf{MIP incumbent} -- the best feasible solution the exact solver
      returns within a budget comparable to the GVNS run time.
\end{itemize}

Because the true optimum is generally not available, the primary quality
measure we report is the \emph{optimality gap of GVNS to the best valid lower
bound},
\[
  gap_{\mathrm{LB}}=\frac{Result_{\text{GVNS}}-\underline{z}}{\underline{z}},
  \qquad
  \underline{z}=\text{best MIP dual bound},
\]
which upper-bounds the true optimality gap of GVNS. We additionally report the
ratio of each reference objective to that of GVNS
($Ratio_{\text{method}}=Result_{\text{method}}/Result_{\text{GVNS}}$), so values
below~$1$ are better than GVNS and above~$1$ are worse.

\paragraph{On the hardness of the exact model.} A first, and perhaps
surprising, observation is that the arc-based MILP is hard even at the smallest
sizes: with the vehicle-count integrality it is a large mixed-integer program
(e.g.\ already $\sim\!6000$ variables at $|V|=10$), and \emph{no} instance --
not even $|V|=10$ -- was solved to \emph{proven} optimality within the budgets
above by either solver. The branch-and-bound \emph{dual bound}, however, is
tight (within a few percent of the incumbent) and, crucially, well above the LP
relaxation; it is therefore the right yardstick for judging GVNS. We also note
that the two solvers behave very differently on this model: their dual bounds
agree closely (mutually cross-validating the lower bound), but HiGHS finds
markedly better incumbents, whereas SCIP fails to find any feasible solution on
several instances with $|V|\!\ge\!45$ and runs out of memory at $|V|=150$. HiGHS
is thus both the stronger exact solver here and the natural in-loop choice.

\begin{figure}[t]
\centering
\includegraphics[width=0.62\textwidth]{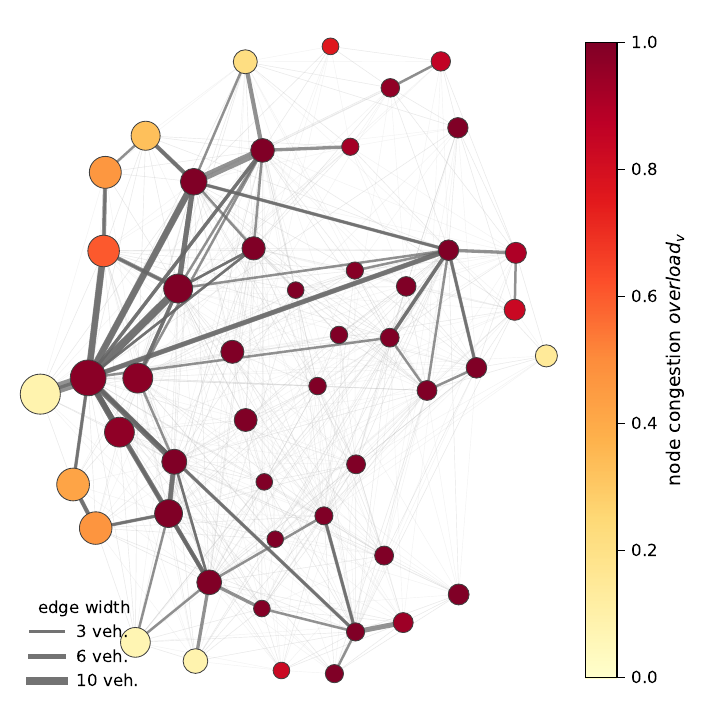}
\caption{A GVNS solution on the Germany real-geometry network ($|V|=49$), drawn
on its geographic layout. Node area encodes the transshipment capacity $W_v$ and
node colour the congestion $overload_v$ (dark $=$ at capacity). All used edges
are drawn faintly; backbone edges carrying at least three vehicles are emphasized
with width proportional to the vehicle count.}
\label{fig:solution}
\end{figure}

\subsection{Results}
\begin{figure}[t]
\centering
\begin{subfigure}{0.49\textwidth}
  \includegraphics[width=\linewidth]{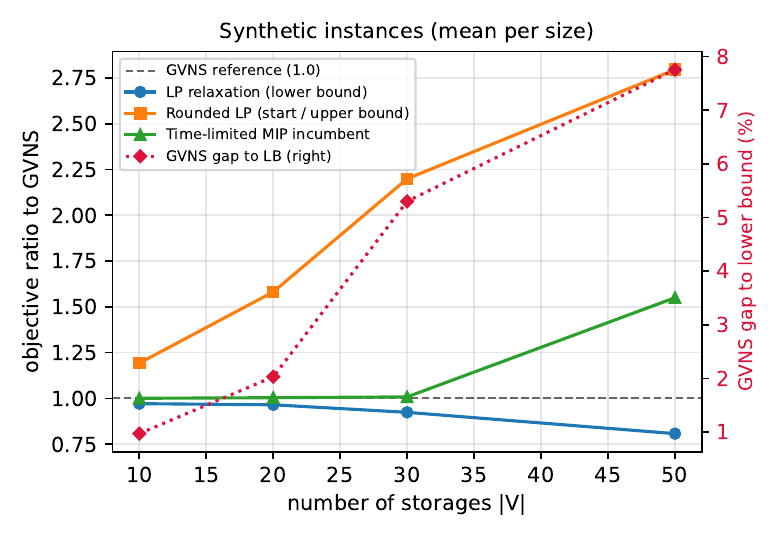}
  \caption{Synthetic instances (mean per size).}
  \label{fig:synth}
\end{subfigure}\hfill
\begin{subfigure}{0.49\textwidth}
  \includegraphics[width=\linewidth]{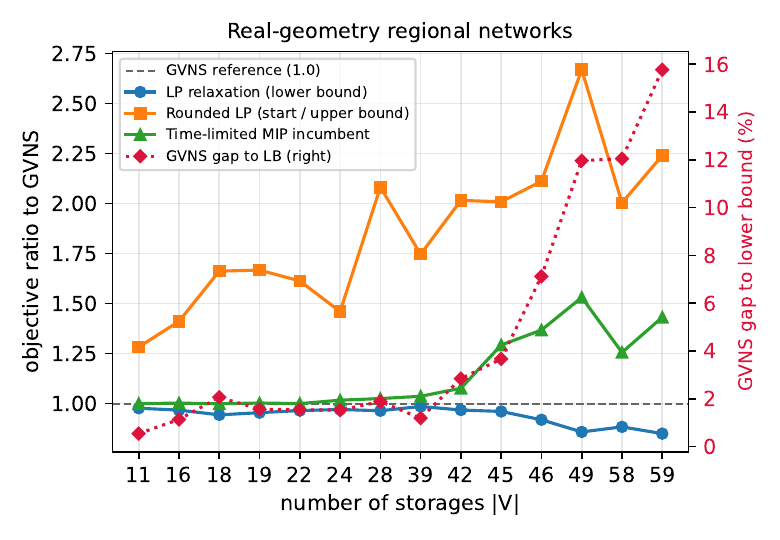}
  \caption{Real-geometry regional networks.}
  \label{fig:real}
\end{subfigure}
\caption{Left axis: objective ratio of each reference method to GVNS versus
network size (values below the GVNS line at $1.0$ are better than GVNS, above
are worse). Right axis (dotted): the GVNS optimality gap to the best valid lower
bound. Only sizes with an exact dual bound are shown.}
\label{fig:ratios}
\end{figure}

Tables~\ref{tab:synth} and~\ref{tab:real} summarise the synthetic and
real-geometry results; Figure~\ref{fig:ratios} plots the optimality gap and the
reference ratios against instance size. Four consistent trends emerge.

\emph{First, GVNS is near-optimal on small instances.} Its gap to the
MIP dual bound is only $0.97\%$ at $|V|=10$ and stays around $1$--$2\%$ up to
$|V|\approx30$, then grows with size (to $\sim\!8\%$ at $|V|=50$), reflecting the
rapidly expanding path space. Since the dual bound is a valid lower bound, this
gap bounds the true distance of GVNS from the optimum rather than comparing it
against itself.

\emph{Second, GVNS is very stable across random seeds.} The across-seed standard
deviation of the objective is below $0.5\%$ of the mean on every instance size
(e.g.\ $0.07\%$ at $|V|=10$, $0.31\%$ at $|V|=50$), so a single run is
representative and the method is robust to its randomisation.

\emph{Third, GVNS matches or beats the exact solver at equal time, increasingly
so with size.} On small instances the time-limited MIP incumbent essentially
ties GVNS ($Ratio_{\text{MIP}}\approx1.00$ for $|V|\le30$); from $|V|=45$ the
exact solver falls behind sharply, returning solutions $1.3$--$1.55\times$ more
expensive than GVNS at $|V|=49$--$59$. In other words, for a comparable
wall-clock budget the exact solver is up to $\sim\!55\%$ worse than GVNS. We are
able to obtain tight dual bounds up to $|V|=59$; for $|V|\gtrsim90$ the arc
model reaches millions of variables and the exact branch-and-bound---both its
tight dual bound and its incumbent---becomes intractable, leaving only the
closed-form but loose LP-relaxation lower bound. Even there GVNS keeps improving
on its rounded-LP start by large factors ($2.2$--$3.0\times$;
Tables~\ref{tab:synth} and~\ref{tab:real}).

\emph{Fourth, the real-geometry networks mirror the synthetic ones}
(Table~\ref{tab:real}): the gap to the dual bound rises from $0.5\%$ (Denmark,
$|V|=11$) to $12$--$16\%$ (Germany $|V|=49$, Spain $|V|=59$), and the exact
incumbent degrades from parity to $1.3$--$1.5\times$ GVNS over the same range.
This corroborates both the fidelity of the generator and the robustness of GVNS.

\begin{table}[t]
\centering
\caption{Synthetic results aggregated by network size. $gap_{\mathrm{LB}}$ is the
GVNS optimality gap to the best MIP dual bound (a valid lower bound);
$Ratio_\bullet$ are objective ratios to GVNS. ``$\dagger$'' = exact model
intractable within the budget (memory/time). $n$ = number of instances.}
\label{tab:synth}
\small
\begin{tabular}{ccccccc}
\toprule
$|V|$ & $n$ & $gap_{\mathrm{LB}}$ & Ratio LP-relax. & Ratio Rounded LP & Ratio MIP-inc. & GVNS std\\
\midrule
10  & 10 & 0.97\% & 0.971 & 1.192 & 1.000 & 0.07\%\\
20  & 7  & 2.03\% & 0.965 & 1.581 & 1.004 & 0.14\%\\
30  & 7  & 5.30\% & 0.924 & 2.199 & 1.008 & 0.49\%\\
50  & 5  & 7.75\% & 0.807 & 2.795 & 1.550 & 0.31\%\\
100 & 3  & \dag & 0.541 & 2.186 & \dag & 0.17\%\\
150 & 3  & \dag & 0.372 & 2.197 & \dag & 0.26\%\\
\bottomrule
\end{tabular}
\\[2pt]
{\footnotesize $\dagger$ at these sizes the arc model has millions of variables,
so the exact branch-and-bound---its tight dual bound and its incumbent---is
intractable; the LP-relaxation lower bound and the rounded-LP upper bound remain
available in closed form (Sections~\ref{sec:relax} and~\ref{sec:gvns}).}
\end{table}

\begin{table}[t]
\centering
\caption{Results on the real-geometry regional networks. Sizes are the storage
counts of the generated instances.}
\label{tab:real}
\small
\begin{tabular}{lcccccc}
\toprule
Network & $|V|$ & $gap_{\mathrm{LB}}$ & Ratio LP-relax. & Ratio Rounded LP & Ratio MIP-inc. & GVNS $t$ (s)\\
\midrule
Denmark        & 11 & 0.5\%  & 0.977 & 1.283 & 1.000 & 523\\
Belgium        & 16 & 1.1\%  & 0.968 & 1.411 & 1.003 & 1279\\
Netherlands    & 18 & 2.1\%  & 0.944 & 1.663 & 0.999 & 2191\\
Hungary        & 19 & 1.6\%  & 0.955 & 1.668 & 1.003 & 975\\
Ohio USA       & 22 & 1.5\%  & 0.965 & 1.613 & 1.000 & 1535\\
Georgia USA    & 24 & 1.5\%  & 0.971 & 1.460 & 1.017 & 1823\\
Illinois USA   & 28 & 1.9\%  & 0.964 & 2.082 & 1.025 & 2605\\
California USA & 39 & 1.2\%  & 0.986 & 1.747 & 1.037 & 2409\\
United Kingdom & 42 & 2.8\%  & 0.968 & 2.016 & 1.077 & 3877\\
Sweden         & 45 & 3.7\%  & 0.961 & 2.009 & 1.293 & 2818\\
Italy          & 46 & 7.1\%  & 0.919 & 2.112 & 1.368 & 2669\\
Germany        & 49 & 12.0\% & 0.858 & 2.670 & 1.531 & 3292\\
Texas USA      & 58 & 12.0\% & 0.884 & 2.004 & 1.256 & 4168\\
Spain          & 59 & 15.8\% & 0.850 & 2.241 & 1.432 & 4145\\
Europe Russia  & 90 & \dag & 0.763 & 3.017 & \dag & 4784\\
\bottomrule
\end{tabular}
\\[2pt]
{\footnotesize $\dagger$ exact branch-and-bound (tight dual bound and incumbent)
intractable at this size; the LP-relaxation and rounded-LP bounds remain
available in closed form.}
\end{table}

\subsection{Neighborhood ablation}
\label{sec:ablation}

\begin{table}[t]
\centering
\caption{Neighborhood ablation: mean objective increase (best of three seeds)
when a single neighborhood is disabled, relative to the full GVNS.}
\label{tab:ablation}
\small
\begin{tabular}{lcc}
\toprule
Disabled neighborhood & Synthetic & Real\\
\midrule
$N_1$ (reload-aware shortest paths) & $+1.5\%$ & $+0.7\%$\\
$N_2$ (partially loaded edges)      & $+7.2\%$ & $+2.5\%$\\
$N_3$ (congestion detours)          & $+0.0\%$ & $+0.8\%$\\
\bottomrule
\end{tabular}
\end{table}
To assess the contribution of each neighborhood we ran a leave-one-out ablation:
GVNS is re-run with a single neighborhood disabled (in both shake and VND) and
compared, over three seeds, against the full method. Table~\ref{tab:ablation}
reports the mean objective increase on a representative subset of synthetic
($|V|=20$--$50$) and real-geometry ($|V|=16$--$49$) instances.

Three effects stand out. $N_2$ (vehicle consolidation over partially loaded
edges) is the dominant operator: removing it costs up to $24.7\%$ on the larger
synthetic instances and $7.2\%$ on Germany, confirming that filling spare vehicle
capacity is the main lever on the objective. $N_1$ contributes a smaller but
consistent improvement. $N_3$ (detours around congested storages) is
\emph{situational}: it is essentially inactive when storages are far from their
limit (the synthetic and small real instances), but becomes valuable exactly when
the transshipment capacity binds---on Germany, where almost every storage operates
at capacity (Figure~\ref{fig:solution}), disabling $N_3$ worsens the objective by
$3.6\%$. Each neighborhood thus earns its place, addressing a different source of
cost.

\section{Conclusion}
We addressed the design of optimal logistic transport flows with nonlinear,
vehicle-induced edge costs and hard storage transshipment capacities. We gave
arc-based and path-based MILP models and their linear relaxation, and proposed a
GVNS that operates in the path space through three problem-specific
neighborhoods: reload-aware shortest paths, partially loaded edges, and detours
around congested storages. 

On $35$ synthetic and $15$ real-geometry instances of up to $150$ storages,
GVNS produced solutions within about one percent of the best dual bound on
small networks (a MIP lower bound far tighter than the LP relaxation) and, for a
comparable budget, up to $55\%$ cheaper than a time-limited exact solver on the
medium ones; it is stable across random seeds, and its optimality gap grows with
size, as expected from the exploding path space. A leave-one-out ablation
confirmed that each neighborhood earns its place: consolidating flow over
partially loaded edges ($N_2$) is the dominant operator, while congestion detours
($N_3$) become valuable specifically on capacity-tight instances. Future work
includes integrating column-generation pricing to feed the path pool and
validating the approach on operational data from real logistic systems.

\subsubsection*{Code and data availability.}
The GVNS implementation, the exact-solver harness and the scripts that reproduce
the experiments of Section~\ref{sec:experiments} are available at
\url{https://github.com/aponom84/Optimal-Logistic-Net}. The benchmark instances
are published as a Zenodo dataset~\cite{benchmark-data} under CC-BY-4.0, and each
is reproducible from its seed and parameters~\cite{benchmark}.

\subsubsection*{Acknowledgements.}
The study was implemented in the framework of the Basic Research Program at
HSE University (HSE-BR-2025-080).
Also, this research was supported in part through computational resources of HPC facilities at HSE University.

\FloatBarrier


\begin{thebibliography}{99}
\small
\bibitem{barnhart1998} C.~Barnhart, E.~L.~Johnson, G.~L.~Nemhauser,
M.~W.~P.~Savelsbergh, P.~H.~Vance. Branch-and-Price: Column Generation for
Solving Huge Integer Programs. \textit{Operations Research}, 46(3):316--329,
1998.
\bibitem{crainic2015} T.~G.~Crainic, M.~Gendreau, J.~M.~Farvolden.
Metaheuristics for Multi-Commodity Flow Problems. \textit{Parallel Computing},
41:1--15, 2015.
\bibitem{masri2015} H.~Masri, S.~Krichen, A.~Guitouni. A multi-start variable
neighborhood search for solving the single path multicommodity flow problem.
\textit{Applied Mathematics and Computation}, 251:132--142, 2015.
\bibitem{amous2017} M.~Amous, S.~Toumi, B.~Jarboui, M.~Eddaly. A variable
neighborhood search algorithm for the capacitated vehicle routing problem.
\textit{Computers \& Operations Research}, 88:234--248, 2017.
\bibitem{hansen2019} P.~Hansen, N.~Mladenovi\'c. Variable Neighborhood Search.
In \textit{Handbook of Metaheuristics}, 3rd ed., Springer, 57--97, 2019.
\bibitem{hansen2017} P.~Hansen, N.~Mladenovi\'c, R.~Todosijevi\'c, S.~Hanafi.
Variable neighborhood search: Basics and variants. \textit{EURO Journal on
Computational Optimization}, 5(3):423--454, 2017.
\bibitem{karakostas2019} P.~Karakostas, A.~Sifaleras, M.~C.~Georgiadis. A
general variable neighborhood search-based solution approach for the
location-inventory routing problem with distribution outsourcing.
\textit{Computers \& Chemical Engineering}, 126:263--279, 2019.
\bibitem{mjirda2014} A.~Mjirda, B.~Jarboui, N.~Mladenovi\'c, C.~Wilbaut,
S.~Hanafi. A general variable neighbourhood search for the multi-product
inventory routing problem. \textit{Int.\ J.\ of Production Research}, 2014.
\bibitem{osm} OpenStreetMap. \url{https://www.openstreetmap.org} (accessed
2026-04-11).
\bibitem{gabaix1999} X.~Gabaix. Zipf's Law for Cities: An Explanation.
\textit{Quarterly Journal of Economics}, 114(3):739--767, 1999.
\bibitem{gastner2006} M.~T.~Gastner, M.~E.~J.~Newman. The spatial structure of
networks. \textit{The European Physical Journal B}, 49(2):247--252, 2006.
\bibitem{barthelemy2011} M.~Barth\'elemy. Spatial networks. \textit{Physics
Reports}, 499(1--3):1--101, 2011.
\bibitem{kirkley2018} A.~Kirkley, H.~Barbosa, M.~Barthelemy, G.~Ghoshal. From
networks to optimal higher-order models of complex systems. \textit{Nature
Physics}, 14:250--258, 2018.
\bibitem{wilson1967} A.~G.~Wilson. A statistical theory of spatial distribution
models. \textit{Transportation Research}, 1(3):253--269, 1967.
\bibitem{highs} Q.~Huangfu, J.~A.~J.~Hall. Parallelizing the dual revised simplex
method. \textit{Mathematical Programming Computation}, 10(1):119--142, 2018.
\bibitem{scip} C.~Hojny, M.~Besan\c{c}on, K.~Bestuzheva, et al. The SCIP
Optimization Suite 10.0. \textit{arXiv preprint} arXiv:2511.18580, 2025.
\bibitem{charisma} Kostenetskiy, P., Kozyrev, V., Chulkevich, R., Salikova, M.,
Mishenin, R.: High-Performance Computing at HSE University. In: Sokolinsky, L.,
Zymbler, M., Voevodin, V., Dongarra, J. (eds.) Parallel Computational
Technologies, pp.~17--29. Springer Nature Switzerland, Cham (2026)
\bibitem{benchmark} V.~Kovaleva, A.~Ponomarenko. A Benchmark Generator of
Realistic Logistic Transport Networks with Nonlinear Edge Costs and
Transshipment Capacities. \textit{arXiv preprint} arXiv:2606.20653, 2026.
\bibitem{benchmark-data} V.~Kovaleva, A.~Ponomarenko. A Benchmark of Realistic
Logistic Transport Networks with Nonlinear Edge Costs and Transshipment
Capacities (dataset). Zenodo, 2026. \url{https://doi.org/10.5281/zenodo.20372456}
\bibitem{batsyn2014ttrp} Batsyn, M., Ponomarenko, A.: Heuristic for a real-life
truck and trailer routing problem. Procedia Computer Science \textbf{31},
778--792 (2014)
\bibitem{batsyn2015heuristic} Batsyn, M., Ponomarenko, A.: Heuristic for
site-dependent truck and trailer routing problem with soft and hard time windows
and split deliveries. In: Pardalos, P., Pavone, M., Farinella, G.M., Cutello, V.
(eds.) Machine Learning, Optimization, and Big Data (MOD 2015). LNCS, vol.~9432,
pp.~65--79. Springer, Cham (2015)
\end{thebibliography}
\end{document}